# A SIMPLE SMOOTH BACKFITTING METHOD FOR ADDITIVE MODELS

By Enno Mammen[1] and Byeong U. Park[2]

*University of Mannheim and Seoul National University*

In this paper a new smooth backfitting estimate is proposed for additive regression models. The estimate has the simple structure of Nadaraya–Watson smooth backfitting but at the same time achieves the oracle property of local linear smooth backfitting. Each component is estimated with the same asymptotic accuracy as if the other components were known.

**1. Introduction.** In additive models it is assumed that the influence of different covariates enters separately into the regression model and that the regression function can be modeled as the sum of the single influences. This is often a plausible assumption. It circumvents fitting of high-dimensional curves and for this reason it avoids the so-called curse of dimensionality. On the other hand, it is a very flexible model that also allows good approximations for more complex structures. Furthermore, the low-dimensional curves fitted in the additive model can be easily visualized in plots. This allows a good data-analytic interpretation of the qualitative influence of single covariates.

In this paper we propose a new backfitting estimate for additive regression models. The estimate is a modification of the smooth backfitting estimate of Mammen, Linton and Nielsen [9]. Their versions of smooth backfitting have been introduced for Nadaraya–Watson smoothing and for local linear smoothing. Smooth backfitting based on Nadaraya–Watson smoothing has the advantage of being easily implemented and of having rather simple intuitive interpretations. On the other hand, local linear smooth backfitting

Received December 2003; revised November 2005.

[1]Supported by DFG Project MA 1026/6-2 of the Deutsche Forschungsgemeinschaft.

[2]Supported by KOSEF, the Korean–German Cooperative Science Program F01-2004-000-10312-0. The paper was written during a visit of B. U. Park to Heidelberg, financed by KOSEF and DFG Project 446KOR113/155/0-1.

leads to more complicated technical implementations. The backfitting formula has no easy interpretation. But, the local linear smooth backfitting estimate has very nice asymptotic properties. It achieves the asymptotic oracle bounds. The local linear smooth backfitting estimate of an additive component has the same asymptotic bias and variance as a theoretical local linear estimate that uses knowledge of the other components. In this paper we introduce a smooth backfitting estimate that has the simple structure of a Nadaraya–Watson estimate but at the same time has the asymptotic oracle property of local linear smoothing.

Several approaches have been proposed for fitting additive models: the classical backfitting procedure by Buja, Hastie and Tibshirani [1], the marginal integration method of Linton and Nielsen [8] and Tjøstheim and Auestad [16], the smooth backfitting estimate of Mammen, Linton and Nielsen [9], the local quasi-differencing approach of Christopeit and Hoderlein [2], and the two-step procedures of Horowitz, Klemelä and Mammen [5]. All estimates require several estimation steps.

The marginal integration estimate makes use of a full-dimensional non-parametric regression estimate as a pilot estimate. Each component of the additive model is fitted by marginal integration of the full-dimensional fit, that is, by integrating out all other arguments. Versions of marginal integration have been proposed that achieve oracle bounds [4]. The algorithm is unstable for moderate and large numbers of additive components and calculation of the full-dimensional regression estimate causes problems. On the other hand, backfitting avoids fitting a full-dimensional regression estimate. It is based on an iterative algorithm. In each step only one additive component is updated. All other components are fixed. So, only one-dimensional smoothing is applied. Asymptotic theory for the classical backfitting is complicated by the fact that the estimate is defined as a limit of the iterative backfitting algorithm but no explicit definition is given. Asymptotic theory is available under restrictive conditions on the design densities [13, 14]. In general, the classical backfitting estimates do not achieve the oracle bounds. For practical implementations of the backfitting estimates, see [15].

Smooth backfitting estimates are defined as the minimizers of a smoothed least squares criterion. As backfitting estimates they can be calculated by an iterative backfitting algorithm. Asymptotic analysis becomes simpler because of the explicit definition of the estimate. Furthermore, making use of an approach in [10], the estimate can be interpreted as an orthogonal projection of the data vector onto the space of additive functions. As with the classical backfitting estimates, smooth backfitting does not make use of a full-dimensional estimate and for this reason it avoids the curse of dimensionality. Smooth backfitting also achieves the oracle bounds. This has been shown for smooth backfitting estimates based on local linear fitting (see [9]). For practical implementations of smooth backfitting, see [12] and [11]. Some two-step procedures have been proposed for additive models. Christopeit and Hoderlein [2] use local quasi-differencing in the second step, an idea coming from efficient estimation in semiparametric estimation. Horowitz, Klemelä and Mammen [5] and Horowitz and Mammen [6] develop a general approach that allows oracle efficient estimates for a broad class of smoothing methods. For a related approach, see also [7].

In the original version of local linear smooth backfitting, both the estimated value and the estimated slope of an additive component are updated. This is done by application of a two-dimensional integral operator. This definition leads to lengthy formulas, which makes it hard to implement the method. Furthermore, the understanding of the method and of its asymptotic properties is complicated by the two-dimensional nature of the integral operator. On the other hand, smooth backfitting for Nadaraya–Watson smoothing is rather simple to understand and it can be rather easily implemented. Again, an integral operator is used in the backfitting steps. But now the operator can be easily interpreted as an empirical analogue of a conditional expectation. In this paper we propose a smooth backfitting estimate that inherits the advantages of Nadaraya–Watson and local linear smooth backfitting. As with Nadaraya–Watson smoothing, it is based on one-dimensional updating. This essentially simplifies the interpretation and asymptotic analysis of the algorithm. On the other hand, the new estimate achieves the asymptotic oracle bounds of local linear smooth backfitting. Our numerical study confirms this asymptotic equivalence, and suggests that the new estimate has a slightly better performance.

The paper is organized as follows. In the next section the method is introduced and is shown to be asymptotically equivalent to local linear smooth backfitting under some conditions on the kernel functions of the backfitting operator. Section 3 discusses some numerical properties of the new proposal. The assumptions for our theoretical results and proofs are deferred to Section 4.

## 2. Local linear smooth backfitting.
In this section we introduce our new smooth backfitting method for local linear smoothing. It is based on a modification of smooth backfitting for Nadaraya–Watson smoothing. We briefly recall the definition of Nadaraya–Watson backfitting from Mammen, Linton and Nielsen [9]. We consider an additive model. For $i = 1, \ldots, n$, it is assumed for one-dimensional response variables $Y^1, \ldots, Y^n$ that

$$(2.1) \qquad Y^i = m_0 + m_1(X_1^i) + \cdots + m_d(X_d^i) + \varepsilon^i.$$

Here, $\varepsilon^i$ are error variables, $m_1, \ldots, m_d$ are unknown functions from $\mathbb{R}$ to $\mathbb{R}$ satisfying $Em_j(X_j^i) = 0$, $m_0$ is an unknown constant and $X^i = (X_1^i, \ldots, X_d^i)$ are random design points in $\mathbb{R}^d$. Throughout the paper we make the assumption that $X^1, \ldots, X^n$ are i.i.d. and that $X^i_j$ takes its values in a bounded interval $I_j$. Furthermore, the error variables $\varepsilon^1, \ldots, \varepsilon^n$ are assumed to be i.i.d. mean zero and to be independent of $X^1, \ldots, X^n$. This excludes interesting autoregression models, but it simplifies our asymptotic analysis. We expect that our results can be extended to dependent observations under mixing conditions.

The Nadaraya–Watson smooth backfitting estimate is defined as the minimizer of the smoothed sum of squares

$$(2.2) \quad \sum_{i=1}^{n} \int_I \left[ Y^i - \widehat{m}_0 - \sum_{j=1}^{d} \widehat{m}_j(x_j) \right]^2 \kappa\left( \frac{X^i_1 - x_1}{h_1}, \ldots, \frac{X^i_d - x_d}{h_d} \right) dx,$$

where $\kappa$ is a $d$-variate kernel function and $I = I_1 \times \cdots \times I_d$. The minimization is done under the constraints

$$(2.3) \quad \int_{I_j} \widehat{m}_j(x_j) \widehat{p}_j(x_j) \, dx_j = 0, \qquad j = 1, \ldots, d,$$

where $\widehat{p}_j$ is a marginal kernel density estimate. The minimizer $\widehat{m}_j$ of (2.2) is uniquely defined by the equations (see [9])

$$(2.4) \quad \widehat{m}_j(x_j) = \widetilde{m}_j(x_j) - \sum_{k \neq j} \int_{I_k} \widehat{m}_k(x_k) \widehat{\pi}_{jk}(x_j, x_k) \, dx_k, \qquad j = 1, \ldots, d,$$

where $\widetilde{m}_j$ is a normalized marginal Nadaraya–Watson estimate and $\widehat{\pi}_{jk}$ is a kernel density estimate of the conditional density $p_{jk}/p_j$. Here $p_{jk}$ denotes the marginal density of $(X_j, X_k)$.

In this paper we propose to use other choices of $\widetilde{m}_j$ and $\widehat{\pi}_{jk}$, and define a new estimate by (2.4) with these new choices. Let $\breve{m}_j$ be the marginal local linear estimate. Together with the slope estimate $\breve{m}^*_j$ the local linear estimate is defined as the minimizer of

$$(2.5) \quad \sum_{i=1}^{n} [Y^i - \breve{m}_j(x_j) - \breve{m}^*_j(x_j)(X^i_j - x_j)]^2 K_{h_j}\,(x_j, X^i_j),$$

where $K_{h_j}$ is a boundary corrected univariate kernel function. It is defined as

$$K_{h_j}(u_j, v_j) = [a(u_j, h_j)v_j + b(u_j, h_j)]h_j^{-1} K\left( \frac{v_j - u_j}{h_j} \right),$$

where $K$ is a symmetric convolution kernel (i.e., a probability density function) supported on $[-1, 1]$. The functions $a$ and $b$ are chosen so that

$$(2.6) \quad \int_{I_j} K_{h_j}(u_j, v_j) \, dv_j = 1,$$

$$(2.7) \quad \int_{I_j} (v_j - u_j) K_{h_j}(u_j, v_j) \, dv_j = 0.$$

We also write $K_{h_j}(v_j - u_j)$ for the kernel $h_j^{-1}K[(v_j - u_j)/h_j]$. This kernel should not be confused with $K_{h_j}(u_j, v_j)$. Specifically,

$$(2.8) \quad \begin{aligned} K_{h_j}(u_j, v_j) &= \left[ \frac{\mu_{K,j,2}(u_j) - (h_j^{-1}(v_j - u_j))\mu_{K,j,1}(u_j)}{\mu_{K,j,0}(u_j)\mu_{K,j,2}(u_j) - \mu_{K,j,1}(u_j)^2} \right] \\ &\quad \times h_j^{-1}K\left( \frac{v_j - u_j}{h_j} \right), \end{aligned}$$

where

$$\mu_{K,j,\ell}(u_j) = \int_{I_j} (v_j - u_j)^\ell h_j^{-\ell} K_{h_j}(v_j - u_j)\, dv_j = \int_{I_j(u_j, h_j, +)} t^\ell K(t)\, dt$$

for $I_j(u_j, h_j, +) = \{t : u_j + h_j t \in I_j\}$.

The normalized marginal estimate $\widetilde{m}_j$ is defined as

$$(2.9) \quad \widetilde{m}_j(x_j) = \check{m}_j(x_j) - \left[ \int \widetilde{p}_j(u)\, du \right]^{-1} \int \check{m}_j(u)\widetilde{p}_j(u)\, du$$

for a modified density estimate $\widetilde{p}_j$. The modified kernel density estimate $\widetilde{p}_j$ is defined as

$$\widetilde{p}_j(x_j) = \widehat{p}_j(x_j) - \frac{\widehat{p}_j^*(x_j)^2}{\widehat{p}_j^{**}(x_j)},$$

where $\widehat{p}_j$ is the usual kernel density estimate,

$$\widehat{p}_j^*(x_j) = \frac{1}{n} \sum_{i=1}^n K_{h_j}(x_j, X_j^i)(X_j^i - x_j),$$

$$\widehat{p}_j^{**}(x_j) = \frac{1}{n} \sum_{i=1}^n K_{h_j}(x_j, X_j^i)(X_j^i - x_j)^2.$$

For the definition of $\widehat{\pi}_{jk}$, we consider the two-dimensional kernel density estimate

$$\widetilde{p}_{jk}(x_j, x_k) = \widehat{p}_{jk}(x_j, x_k) - \frac{\widehat{p}_{jk}^*(x_j, x_k)\widehat{p}_j^*(x_j)}{\widehat{p}_j^{**}(x_j)},$$

where

$$\widehat{p}_{jk}(x_j, x_k) = \frac{1}{n} \sum_{i=1}^n K_{h_j}(x_j, X_j^i)L_{h_k}(x_k, X_k^i),$$

$$\widehat{p}_{jk}^*(x_j, x_k) = \frac{1}{n} \sum_{i=1}^n K_{h_j}(x_j, X_j^i)L_{h_k}(x_k, X_k^i)(X_j^i - x_j).$$

The kernel $L_{h_k}$ is defined as

$$L_{h_k}(u_k, v_k) = [c(v_k, h_k)u_k + d(v_k, h_k)]h_k^{-1} L\left(\frac{v_k - u_k}{h_k}\right),$$

where $c$ and $d$ are chosen so that

(2.10)                              $$\int_{I_k} L_{h_k}(u_k, v_k)\, du_k = 1,$$

(2.11)                              $$\int_{I_k} (v_k - u_k) L_{h_k}(u_k, v_k)\, du_k = 0.$$

Specifically,

$$L_{h_k}(u_k, v_k) = \left[\frac{\mu^*_{L,k,2}(v_k) - (h_k^{-1}(v_k - u_k))\mu^*_{L,k,1}(v_k)}{\mu^*_{L,k,0}(v_k)\mu^*_{L,k,2}(v_k) - \mu^*_{L,k,1}(v_k)^2}\right] h_k^{-1} L\left(\frac{v_k - u_k}{h_k}\right),$$

where

$$\mu^*_{L,k,\ell}(v_k) = \int_{I_k} (v_k - u_k)^\ell h_k^{-\ell} L_{h_k}(v_k - u_k)\, du_k = \int_{I_k(v_k, h_k, -)} t^\ell L(t)\, dt$$

for $I_k(v_k, h_k, -) = \{t : v_k - h_k t \in I_k\}$. We use the following convolution kernel $L$:

$$L(u) = 2K_{1/\sqrt{2}}(u) - K_{\sqrt{2}}(u).$$

This kernel satisfies $\int L(u)\, du = 1$, $\int u L(u)\, du = 0$ and $\int u^2 L(u)\, du = -\int u^2 \times K(u)\, du$. Other kernels with these moments will also work. Again, we also write $L_{h_j}(v_j - u_j)$ for the kernel $h_j^{-1} L[(v_j - u_j)/h_j]$. Note that the definition of $L_{h_k}$ differs from that of $K_{h_j}$. The difference comes from integration with respect to different arguments in the moment equations. Note also that the moment condition (2.10) is required on their kernels $K_{h_k}$ (as well as $K_{h_j}$) for the local linear smooth backfitting estimate proposed by Mammen, Linton and Nielsen [9]. The additional condition (2.11) on the first moment is needed here to mimic local linear smooth backfitting with a Nadaraya–Watson-type estimate.

We now define $\widehat{\pi}_{jk}$ as

(2.12)                     $$\widehat{\pi}_{jk}(x_j, x_k) = \frac{\widetilde{p}_{jk}(x_j, x_k)}{\widetilde{p}_j(x_j)} - \frac{\int \widetilde{p}_{jk}(u, x_k)\, du}{\int \widetilde{p}_j(u)\, du}.$$

Our main result states that the estimate $\widehat{m}_j$ is asymptotically equivalent to local linear smooth backfitting estimates. We will give motivation for the choice of $\widehat{\pi}_{jk}$ at the end of this section.

THEOREM 2.1. *Under Assumptions* (A1)–(A5) *stated in Section* 4, *we get the following expansions for the estimate* $\widehat{m}_j$ *defined by* (2.4) *with* $\widetilde{m}_j$ *at* (2.5)–(2.9) *and* $\widehat{\pi}_{jk}$ *at* (2.12):

$$
\begin{aligned}
(2.13) \quad \widehat{m}_j(x_j) &= m_j(x_j) + h_j^2 \Big[ \frac{1}{2} C_{K,j,2}(x_j) m_j''(x_j) + \Delta_j \Big] \\
&\quad + \Big[ \frac{1}{n} \sum_{i=1}^n K_{h_j}(x_j, X_j^i) \Big]^{-1} \frac{1}{n} \sum_{i=1}^n K_{h_j}(x_j, X_j^i) \varepsilon^i + o_p(n^{-2/5})
\end{aligned}
$$

*uniformly for* $x_j \in I_j$, *where* $C_{K,j,\ell}(x_j) = \int_{I_j} (v_j - x_j)^\ell h_j^{-\ell} K_{h_j}(x_j, v_j) \, dv_j$,

$$
\begin{aligned}
\Delta_j &= -\frac{1}{2} C_K \Big[ \int m_j(u_j) p_j''(u_j) \, du_j \\
&\quad - 2 \int m_j(u_j) \frac{p_j'(u_j)^2}{p_j(u_j)} \, du_j + \int_{I_j} m_j''(u_j) p_j(u_j) \, du_j \Big]
\end{aligned}
$$

*and* $C_K = \int u^2 K(u) \, du$.

We point out that $C_{K,j,\ell}(x_j)$ in the theorem is different from $\mu_{K,j,\ell}(x_j)$ defined earlier. In fact, for $K_{h_j}$ satisfying (2.6) and (2.7) it follows that

$$
C_{K,j,\ell}(x_j) = \frac{\mu_{K,j,2}(x_j) \mu_{K,j,\ell}(x_j) - \mu_{K,j,1}(x_j) \mu_{K,j,\ell+1}(x_j)}{\mu_{K,j,0}(x_j) \mu_{K,j,2}(x_j) - \mu_{K,j,1}(x_j)^2}.
$$

We compare the estimate $\widehat{m}_j$ with the local linear smooth backfitting estimate, $\widehat{m}_{j,\mathrm{SB}}$, studied by Mammen, Linton and Nielsen [9]. There are differences at the boundary and in the interior of $I_j$. For $x_j$ in the interior $I_j^- = \{u \in I_j : u + h_j \in I_j, \ u - h_j \in I_j\}$ one gets $C_{K,j,2}(x_j) = C_K$. Thus the expansion of $\widehat{m}_j$ becomes

$$
\begin{aligned}
(2.14) \quad \widehat{m}_j(x_j) &= m_j(x_j) + \frac{1}{2} C_K h_j^2 m_j''(x_j) + h_j^2 \Delta_j \\
&\quad + \Big[ \frac{1}{n} \sum_{i=1}^n K_{h_j}(x_j, X_j^i) \Big]^{-1} \frac{1}{n} \sum_{i=1}^n K_{h_j}(x_j, X_j^i) \varepsilon^i + o_p(n^{-2/5}).
\end{aligned}
$$

For $x_j \in I_j^-$ this expansion differs from the stochastic expansion of $\widehat{m}_{j,\mathrm{SB}}$ only by the *constant* term $h_j^2 \Delta_j$; see [9] and [11]. This additive term comes from the norming $\int_{I_j} \widehat{m}_j(u_j) \widetilde{p}_j(u_j) \, du_j = 0$. This can be easily verified by observing that

$$
\begin{aligned}
\int_{I_j} m_j(u_j) &\widetilde{p}_j(u_j) \, du_j \\
&= \frac{1}{2} C_K h_j^2 \Big[ \int m_j(u_j) p_j''(u_j) \, du_j - 2 \int m_j(u_j) \frac{p_j'(u_j)^2}{p_j(u_j)} \, du_j \Big] + o_p(n^{-2/5}).
\end{aligned}
$$

One could use other normings for estimation of $m_j$. We briefly discuss two other normings,

$$(2.15) \qquad \widehat{m}_{j,+}(x_j) = \widehat{m}_j(x_j) - \int_{I_j} \widehat{m}_j(u_j) \widehat{p}_j(u_j) \, du_j,$$

$$(2.16) \qquad \widehat{m}_{j,++}(x_j) = \widehat{m}_j(x_j) - \frac{1}{n} \sum_{i=1}^n \widehat{m}_j(X_j^i).$$

For these two modified estimates the following expansions hold.

COROLLARY 2.1.   *Under the assumptions of Theorem* 2.1, *the expansion* (2.13) *applies for the estimates* $\widehat{m}_{j,+}(x_j)$ *and* $\widehat{m}_{j,++}(x_j)$ *defined at* (2.15) *and* (2.16), *respectively, with* $\Delta_j$ *replaced by*

$$\Delta_{j,+} = -\tfrac{1}{2} C_K \left[ \int m_j(u_j) p_j''(u_j) \, du_j + \int_{I_j} m_j''(u_j) p_j(u_j) \, du_j \right]$$

*for* $\widehat{m}_{j,+}$ *and by*

$$\Delta_{j,++} = -\tfrac{1}{2} C_K \left[ \int_{I_j} m_j''(u_j) p_j(u_j) \, du_j \right]$$

*for* $\widehat{m}_{j,++}$.

For the local linear smooth backfitting estimate $\widehat{m}_{j,\mathrm{SB}}$, we get the expansion (2.14) with $\Delta_{j,\mathrm{SB}} = 0$ for $x_j$ in the interior $I_j^-$; see [9] and [11]. There a different norming was used for a combination of the smooth backfitting estimate of $m_j$ and its derivative; see (3.4) in [11]. The norming of $\widehat{m}_{j,++}$ is chosen so that the mean integrated squared error $\int_{I_j} [\widehat{m}_{j,++}(x_j) - m(x_j)]^2 p(x_j) \, dx_j$ is asymptotically minimized. Note that

$$\int_{I_j} [\widehat{m}_{j,++}(x_j) - m_j(x_j)] p(x_j) \, dx_j = \int_{I_j} \widehat{m}_{j,++}(x_j) p(x_j) \, dx_j = o_P(n^{-2/5}).$$

Furthermore, our estimate $\widehat{m}_j$ differs from the local linear smooth backfitting estimate $\widehat{m}_{j,\mathrm{SB}}$ on the boundary $I_j \backslash I_j^-$. The estimates have slightly different asymptotic biases on the boundary. The difference is due to the fact that they use different boundary corrected kernels. Recall that the local linear estimate in the univariate nonparametric regression with a conventional kernel $K$, without boundary modification, has the asymptotic bias

$$\frac{1}{2} m''(x) \frac{\mu_{K,2}(x)^2 - \mu_{K,1}(x) \mu_{K,3}(x)}{\mu_{K,0}(x) \mu_{K,2}(x) - \mu_{K,1}(x)^2} h^2;$$

see [3]. Here $m$ is the nonparametric regression function, $h$ is the bandwidth, $\mu_{K,\ell}(x) = \int_I (u - x)^\ell h^{-\ell} K_h(u - x) \, du$ for $\ell \geq 0$ and $I$ is the support of the covariate. A similar bias expansion holds for the local linear smooth backfitting estimate $\hat{m}_{j,\mathrm{SB}}(x_j)$. Recall that in the construction of $\hat{m}_{j,\mathrm{SB}}$, boundary corrected kernels $K^*_{h_k}$ that satisfy $\int_{I_k} K^*_{h_k}(x_k, v_k)\, dx_k = 1$ for all $k$ (including $j$) are used. Note that this moment condition is different from (2.6) but is the same as (2.10) that we require on $L$. By an extension of the arguments given in [9] and [11], one gets for the bias of $\hat{m}_{j,\mathrm{SB}}(x_j)$ the expansion

$$\frac{1}{2} m''_j(x_j) \frac{C_{K^*,j,2}(x_j)^2 - C_{K^*,j,1}(x_j) C_{K^*,j,3}(x_j)}{C_{K^*,j,0}(x_j) C_{K^*,j,2}(x_j) - C_{K^*,j,1}(x_j)^2} h_j^2,$$

where $C_{K^*,j,\ell}$ is defined in the same way as $C_{K,j,\ell}$ but with $K_{h_j}$ being replaced by $K^*_{h_j}$. The bias expansion of our estimate is simplified since $C_{K,j,0} = 1$ and $C_{K,j,1} = 0$ from (2.6) and (2.7), respectively.

The asymptotic variances of our estimate $\hat{m}_j$ and the local linear smooth backfitting estimate $\hat{m}_{j,\mathrm{SB}}$ are also slightly different on the boundary. They are identical in the interior of $I_j$, however. The difference on the boundary is also due to the use of different kernels as is discussed above.

We now give motivation for our choice of $\hat{\pi}_{jk}$ when $d = 2$. We give some heuristic arguments why our proposal is a second-order modification of local linear smooth backfitting. We restrict the discussion to points in the interior of $I_j$ and for simplicity we neglect boundary effects. For this reason in the heuristics we use convolution kernels that are not corrected at the boundary. The local linear smooth backfitting estimate of Mammen, Linton and Nielsen [9] is defined as the minimizer of

$$(2.17) \quad \begin{aligned} \sum_{i=1}^{n} \int [Y^i - \hat{m}_0 - \hat{m}_1(x_1) - \hat{m}^*_1(x_1)(X^i_1 - x_1) - \hat{m}_2(x_2) \\ - \hat{m}^*_2(x_2)(X^i_2 - x_2)]^2 K_{h_1}(X^i_1 - x_1) K_{h_2}(X^i_2 - x_2)\, dx_1\, dx_2. \end{aligned}$$

Here $\hat{m}_1$ and $\hat{m}_2$ are the estimates of the additive components $m_1$ and $m_2$, respectively, and $\hat{m}^*_1$ and $\hat{m}^*_2$ are the estimates of the slopes of $m_1$ and $m_2$. Minimization of (2.17) with respect to $\hat{m}_1(x_1)$ and $\hat{m}^*_1(x_1)$ for fixed $x_1$ and for fixed functions $\hat{m}_2(\cdot), \hat{m}^*_2(\cdot)$ leads to

$$(2.18) \quad \begin{aligned} 0 = \sum_{i=1}^{n} \int [Y^i - \hat{m}_0 - \hat{m}_1(x_1) \\ - \hat{m}^*_1(x_1)(X^i_1 - x_1) - \hat{m}_2(x_2) - \hat{m}^*_2(x_2)(X^i_2 - x_2)] \\ \times \begin{pmatrix} 1 \\ X^i_1 - x_1 \end{pmatrix} K_{h_1}(X^i_1 - x_1) K_{h_2}(X^i_2 - x_2)\, dx_2. \end{aligned}$$

This equation is used in the smooth backfitting algorithm for updating $\hat{m}_1$ and $\hat{m}^*_1$. We modify this equation so that the slope estimates $\hat{m}^*_1$ and $\hat{m}^*_2$ do not enter the updating equation and thus the algorithm only keeps track of the values of $\widehat{m}_1$ and $\widehat{m}_2$.

We first discuss how $\widehat{m}_2^*$ can be dropped. The basic idea is to replace equation (2.18) by

$$
\begin{aligned}
0 = \sum_{i=1}^n \int & [Y^i - \widehat{m}_0 - \widehat{m}_1(x_1) - \widehat{m}_1^*(x_1)(X_1^i - x_1) - \widehat{m}_2(x_2)] \\
& \times \binom{1}{X_1^i - x_1} K_{h_1}(X_1^i - x_1) L_{h_2}(X_2^i - x_2) \, dx_2.
\end{aligned}
\tag{2.19}
$$

Here, $L_{h_2}$ is a kernel such that the right-hand sides of (2.18) and (2.19) are asymptotically equivalent. This can be accomplished by choosing $L_{h_2}$ so that

$$
\int [\widehat{m}_2(x_2) + \widehat{m}_2^*(x_2)(X_2^i - x_2)] K_{h_2}(X_2^i - x_2) \, dx_2 \simeq \int \widehat{m}_2(x_2) L_{h_2}(X_2^i - x_2) \, dx_2.
$$

This is done if we choose $L$ to satisfy $\int L(u) \, du = 1$, $\int u L(u) \, du = 0$ and $\int u^2 L(u) \, du = - \int u^2 K(u) \, du$ since

$$
\begin{aligned}
\widehat{m}_2(x_2) &\simeq \widehat{m}_2(X_2^i) - \widehat{m}_2'(X_2^i)(X_2^i - x_2) + \tfrac{1}{2}\widehat{m}_2''(X_2^i)(X_2^i - x_2)^2, \\
\widehat{m}_2^*(x_2) &\simeq \widehat{m}_2'(X_2^i) - \widehat{m}_2''(X_2^i)(X_2^i - x_2).
\end{aligned}
$$

It remains to modify (2.19) further so that $\widehat{m}_1^*$ does not appear. This can be easily achieved by solving (2.19) with respect to $\widehat{m}_1$. It gives

$$
\begin{aligned}
\widehat{m}_1(x_1) = & \left[ \sum_{i=1}^n (X_1^i - x_1)^2 K_{h_1}(X_1^i - x_1) \sum_{i=1}^n Z^i K_{h_1}(X_1^i - x_1) \right. \\
& \left. - \sum_{i=1}^n (X_1^i - x_1) K_{h_1}(X_1^i - x_1) \sum_{i=1}^n Z^i (X_1^i - x_1) K_{h_1}(X_1^i - x_1) \right] \\
& \times \left[ \sum_{i=1}^n (X_1^i - x_1)^2 K_{h_1}(X_1^i - x_1) \sum_{i=1}^n K_{h_1}(X_1^i - x_1) \right. \\
& \left. - \left( \sum_{i=1}^n (X_1^i - x_1) K_{h_1}(X_1^i - x_1) \right)^2 \right]^{-1} - \widehat{m}_0
\end{aligned}
$$

with $Z^i = Y^i - \int \widehat{m}_2(x_2) L_{h_2}(X_2^i - x_2) \, dx_2$. This is equivalent to

$$
\widehat{m}_1(x_1) = \breve{m}_1(x_1) - \widehat{m}_0 - \int \widehat{m}_2(x_2) \frac{\widetilde{p}_{12}(x_1, x_2)}{\widetilde{p}_1(x_1)} \, dx_2,
$$

which implies

$$
\widehat{m}_1(x_1) = \widetilde{m}_1(x_1) - \int \widehat{m}_2(x_2) \widehat{\pi}_{12}(x_1, x_2) \, dx_2.
$$

The above argument is valid for $x_j \in I_j^-$. For the boundary area $I_j \setminus I_j^-$, it continues to hold if one uses the boundary corrected kernel $L_{h_2}(x_2, X_2^i)$ instead of $L_{h_2}(X_2^i - x_2)$ and $K_{h_1}(x_1, X_1^i)$ instead of $K_{h_1}(X_1^i - x_1)$.

**3. Numerical properties.** In this section we compare some numerical properties of the new and the local linear smooth backfitting estimates. For this, we drew 500 datasets $(X^i, Y^i)$, $i = 1, \ldots, n$, with $n = 100$ and $400$ from the model

$$\text{(M1)} \qquad Y^i = m_1(X_1^i) + m_2(X_2^i) + m_3(X_3^i) + \varepsilon^i,$$

where $m_1(x_1) = x_1^2$, $m_2(x_2) = x_2^3$, $m_3(x_3) = -x_3^4$ and $\varepsilon^i$ are distributed as $N(0, 0.01)$. The covariates were generated from truncated three-dimensional normal distributions with marginals $N(0.5, 0.5)$ and correlations $\rho_{12} = \rho_{13} = \rho_{23} = \rho$, where $\rho_{ij}$ denotes the correlation between $X_i$ and $X_j$. The truncation was done for the covariates to have the compact support $[0, 1]^3$. To be specific, a random variate generated from one of the three-dimensional normal distributions was discarded if one of the covariates fell outside the interval $[0, 1]$. The correlation levels used were $\rho = 0$ and $0.5$. The kernel that we used for the backfitting algorithm was the biweight kernel $K(u) = (15/16)(1 - u^2)^2 I_{[-1,1]}(u)$. For the local linear smooth backfitting estimate, we used $K_{h_j}$ that satisfy $\int K_{h_j}(u, v) \, du = 1$ for all $j$, but neither (2.6) nor (2.7). For a fair comparison, we used for the new estimate the conventional kernels $K_{h_j}(v - u)$ instead of $K_{h_j}(u, v)$ given in (2.8). Also, both the new and the local linear smooth backfitting estimates were recentered according to the formula (2.16).

Figures 1 and 2 and Table 1 summarize the results. The target functions are $m_j - Em_j(X_j)$ rather than $m_j$ since $Em_j(X_j) \neq 0$. Figures 1 and 2 depict the bias, the variance and the mean squared error curves of the new and the local linear smooth backfitting estimates, which are based on 500 pseudosamples of size 400. The results for the samples of size 100 are not presented here, but they give the same message as those for the samples of size 400. Table 1 shows the integrated squared biases, integrated variances and integrated mean squared errors. It is observed from Figures 1 and 2 that the bias property of the new estimate $\widehat{m}_j$ is nearly the same as that of the local linear smooth backfitting estimate $\widehat{m}_{j,\mathrm{SB}}$ in the interior and on the boundary. In the interior the variance properties of the two estimates are also nearly the same, while on the boundary the new estimate is seen to be slightly more stable. Because of less variability on the boundary, the new estimate has a slightly improved mean integrated squared error property, as shown in Table 1.

The bandwidths $h_j$ used for these results were chosen as

$$\text{(3.1)} \qquad h_j = n^{-1/5} \left[ E(\varepsilon^i)^2 \int K^2(t) \, dt \right]^{1/5} \left[ C_K^2 \int_0^1 m_j''(u_j)^2 p_j(u_j) \, du_j \right]^{-1/5}.$$

This is the optimal bandwidth for local linear smoothing in univariate non-parametric regression models (i.e., additive models with one additive component) and also the optimal bandwidth for the local linear smooth backfitting estimate $\widehat{m}_{j,\mathrm{SB}}$; see [11] for the latter. In additive models the optimal bandwidth depends on the norming of the estimate. In particular, for the MISE-optimal norming we get the estimate $\widehat{m}_{j,++}(x_j)$ (see the discussion after Theorem 2.1) and an asymptotically optimal bandwidth that is defined as in (3.1) but with $m_j''(u_j)$ replaced by $m_j''(u_j) - \int_{I_j} m_j''(v_j)p(v_j)\,dv_j$. We used the bandwidth as defined in (3.1). In this respect we follow the usual practice in classical nonparametric regression and do not minimize MISE by using estimates of $\int m_j(u)p_j(u)\,du$ that have parametric rate $n^{-1/2}$. Note that in univariate nonparametric regression an estimate $\widehat{m}(x)$ could be improved by the modification $\widehat{m}^*(x) = \widehat{m}(x) - \frac{1}{n}\sum_{i=1}^{n}\widehat{m}(X^i) + n^{-1}\sum_{i=1}^{n}Y^i$. For example, if $\widehat{m}(x)$ is the local linear smoother, then the asymptotic bias of $\widehat{m}^*(x)$ equals $\frac{1}{2}C_K[m''(x) - \int m''(u)p(u)\,du]h^2$, leading to a smaller first-order integrated squared bias. We tried other fixed bandwidths around the optimal bandwidth (3.1), but the lessons were essentially the same.

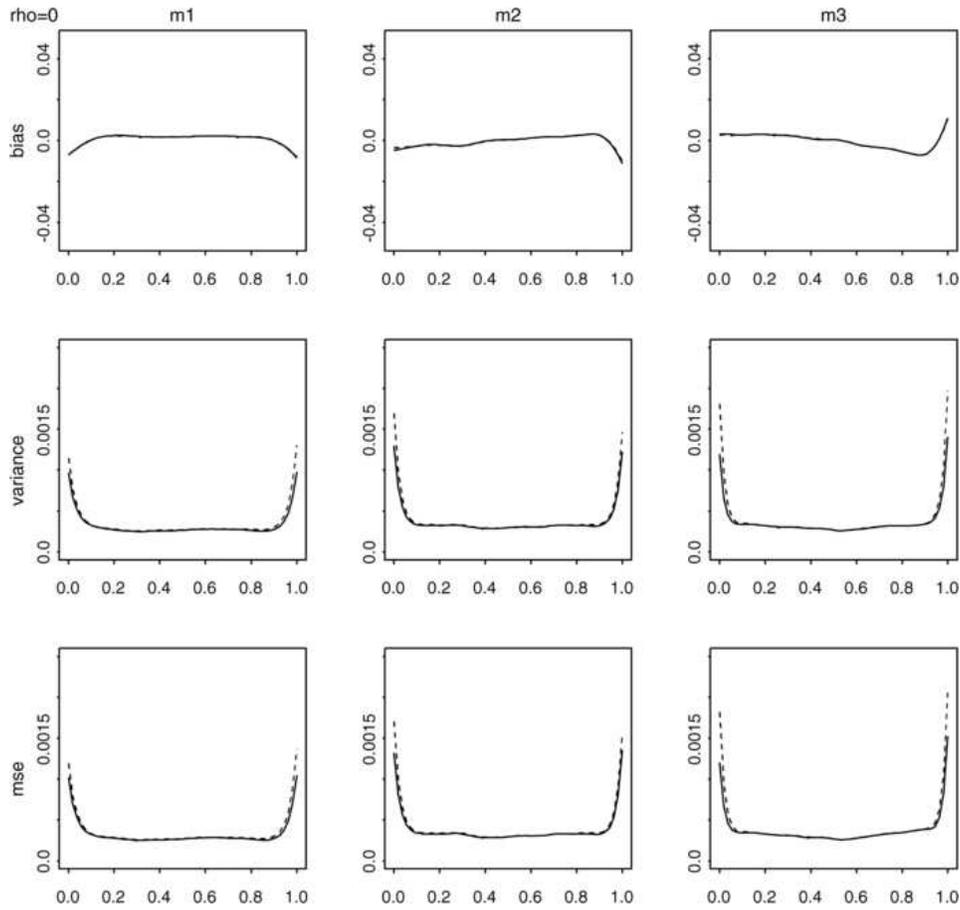

FIG. 1. *Bias, variance and mean squared error curves when $\rho = 0$. The solid curves correspond to the new estimates $\widehat{m}_j$, and the dashed curves are for the local linear smooth backfitting estimates $\widehat{m}_{j,\mathrm{SB}}$. The three rows show the bias, the variance and the mean squared error curves. In each row, the leftmost panel corresponds to $m_1$ and the next two to the right are for $m_2$ and $m_3$. These are based on 500 pseudosamples of size 400.*

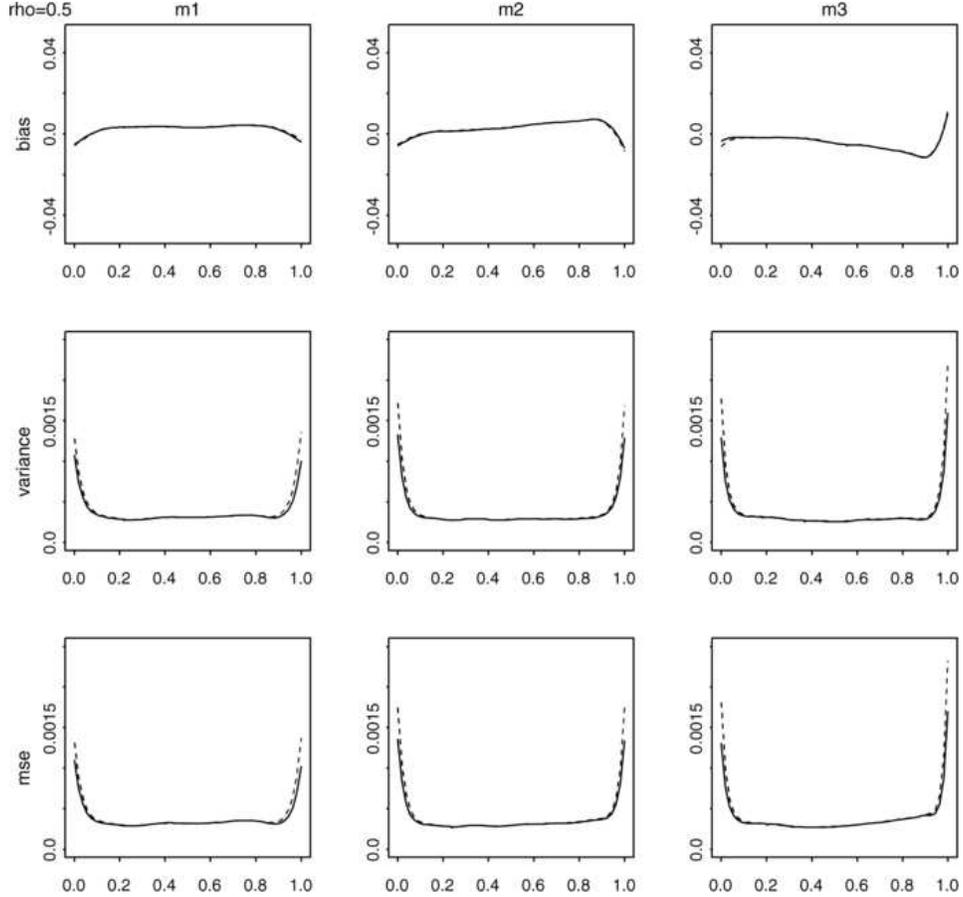

Fig. 2. *Bias, variance and mean squared error curves when $\rho = 0.5$. Line types and arrangement of panels are the same as in Figure* 1. *These are also based on* 500 *pseudosamples of size* 400.

## 4. Assumptions and proofs.

### 4.1. *Assumptions.* Below, we collect the assumptions used in this paper.

Table 1
*Integrated squared bias, integrated variance and integrated mean squared error, multiplied by $10^3$, of the new and the local linear smooth backfitting estimates based on 500 pseudosamples of size 400*

| Corr. level | Target function | Estimate | Integrated sq. bias | Integrated variance | Integrated MSE |
|---|---|---|---|---|---|
| $\rho = 0$ | $m_1$ | $\widehat{m}_{1,\text{SB}}$ | 0.0073 | 0.3479 | 0.3552 |
| | | $\widehat{m}_1$ | 0.0071 | 0.3234 | 0.3305 |
| | $m_2$ | $\widehat{m}_{2,\text{SB}}$ | 0.0070 | 0.4027 | 0.4097 |
| | | $\widehat{m}_2$ | 0.0081 | 0.3768 | 0.3849 |
| | $m_3$ | $\widehat{m}_{3,\text{SB}}$ | 0.0136 | 0.4040 | 0.4176 |
| | | $\widehat{m}_3$ | 0.0138 | 0.3660 | 0.3798 |
| $\rho = 0.5$ | $m_1$ | $\widehat{m}_{1,\text{SB}}$ | 0.0114 | 0.3910 | 0.4024 |
| | | $\widehat{m}_1$ | 0.0114 | 0.3657 | 0.3771 |
| | $m_2$ | $\widehat{m}_{2,\text{SB}}$ | 0.0179 | 0.3928 | 0.4107 |
| | | $\widehat{m}_2$ | 0.0179 | 0.3629 | 0.3808 |
| | $m_3$ | $\widehat{m}_{3,\text{SB}}$ | 0.0334 | 0.3967 | 0.4301 |
| | | $\widehat{m}_3$ | 0.0326 | 0.3601 | 0.3927 |

(A1) The kernel $K$ is bounded and symmetric about zero. It has compact support ($[-1, 1]$, say) and is Lipschitz continuous.

(A2) The $d$-dimensional vector $X^i$ has compact support $I = I_1 \times \cdots \times I_d$ for bounded intervals $I_j$ and its density is bounded from below and from above on $I$.

(A3) $E(\varepsilon^i)^2 < +\infty$.

(A4) The functions $m''_j$, $p''_j$, $D^2_{x_j} p_{jk}(x_j, x_k)$ for $1 \le j$, $k \le d$ exist and are continuous, where $D_{x_j}$ denotes the partial derivative operator with respect to $x_j$ and $D^2_{x_j}$ is the operator of order 2.

(A5) The bandwidths $h_1, \ldots, h_d$ are of order $n^{-1/5}$.

4.2. *Proof of Theorem* 2.1. Define $\eta^i_{kj} = m_k(X^i_k) - E[m_k(X^i_k)|X^i_j]$ and

$$\check{m}^A_j(x_j) = \frac{1}{n} \sum_{i=1}^{n} \frac{K_{h_j}(x_j, X^i_j)}{\widehat{p}_j(x_j)} \varepsilon^i,$$

$$\check{m}^B_j(x_j) = m_0 + m_j(x_j) + \sum_{k \neq j} \int_{I_k} \frac{p_{jk}(x_j, x_k)}{p_j(x_j)} m_k(x_k)\,dx_k$$

$$(4.1) \qquad\qquad + \frac{1}{2} C_{K,j,2}(x_j) h^2_j$$

$$\times \left[ m''_j(x_j) + \sum_{k \neq j} D^2_{x_j} \int \frac{p_{jk}(x_j, x_k)}{p_j(x_j)} m_k(x_k)\,dx_k \right]$$

$$+ \frac{1}{n} \sum_{i=1}^{n} \frac{K_{h_j}(x_j, X_j^i)}{\widehat{p}_j(x_j)} \sum_{k \neq j} \eta_{kj}^i.$$

For the local linear estimate $\check{m}_j$, the following expansion holds:

$$(4.2) \qquad \check{m}_j(x_j) = \check{m}_j^A(x_j) + \check{m}_j^B(x_j) + o_p(n^{-2/5})$$

uniformly for $x_j \in I_j$. These expansions follow by standard asymptotic smoothing theory. Define now

$$\overline{m}_j^B(x_j) = m_j(x_j) + \tfrac{1}{2} C_{K,j,2}(x_j) h_j^2 m_j''(x_j),$$

$$\overline{m}_j^A(x_j) = \check{m}_j^A(x_j).$$

We will show that for $S = A, B$

$$(4.3) \quad \overline{m}_j^S(x_j) = \check{m}_j^S(x_j) - \sum_{k \neq j} \int_{I_k} \overline{m}_k^S(x_k) \frac{\widetilde{p}_{jk}(x_j, x_k)}{\widetilde{p}_j(x_j)} \, dx_k + o_p(n^{-2/5})$$

uniformly for $x_j \in I_j$. Below we argue that (4.3) implies the statement of Theorem 2.1. The proof of (4.3) will be given afterwards.

We apply Theorems 2 and 3 in [9]. We will do this with our $\widetilde{p}_{jk}$, $\widetilde{p}_j$, $\widehat{m}_j^S$, $\check{m}_j^S$, respectively, taking the roles of their $\widehat{p}_{jk}$, $\widehat{p}_j$, $\widetilde{m}_j^S$, $\widehat{m}_j^S$. It is easy to verify the conditions of these theorems. From their Theorem 2 with $S_j = I_j$ and $\Delta_n = n^{-2/5}$ together with our (4.3) we get

$$(4.4) \quad \widehat{m}_j^A(x_j) = \check{m}_j^A(x_j) - \left[ \int \widetilde{p}_j(u) \, du \right]^{-1} \int \check{m}_j^A(u) \widetilde{p}_j(u) \, du + o_p(n^{-2/5})$$

uniformly for $x_j \in I_j$. Here for $S = A, B$ the random function $\widehat{m}_j^S$ is defined by

$$\widehat{m}_j^S(x_j) = \check{m}_j^S(x_j) - \left[ \int \widetilde{p}_j(u) \, du \right]^{-1} \int \check{m}_j^S(u) \widetilde{p}_j(u) \, du$$

$$- \sum_{k \neq j} \int_{I_k} \widehat{m}_k^S(x_k) \widehat{\pi}_{jk}(x_j, x_k) \, dx_k,$$

$$\int \widehat{m}_j^S(u) \widetilde{p}_j(u) \, du = 0.$$

It is easy to check that the second term on the right-hand side of (4.4) is of order $o_p(n^{-2/5})$. Therefore we have

$$(4.5) \qquad \widehat{m}_j^A(x_j) = \check{m}_j^A(x_j) + o_p(n^{-2/5}).$$

Note that

$$(4.6) \qquad \widehat{m}_j(x_j) = \widehat{m}_j^A(x_j) + \widehat{m}_j^B(x_j).$$

We now apply Theorem 3 in [9] with $\alpha_{n,j}(x_j) = \overline{m}_j^B(x_j)$, $\beta(x) \equiv 0$, $\widehat{\mu}_{n,0} = 0$, $\alpha_{n,0} = 0$, $S_j = I_j$ and $\Delta_n = n^{-2/5}$. This gives

$$(4.7) \quad \widehat{m}_j^B(x_j) = \overline{m}_j^B(x_j) - \left[ \int \widetilde{p}_j(u) \, du \right]^{-1} \int \overline{m}_j^B(u) \widetilde{p}_j(u) \, du + o_p(n^{-2/5})$$

uniformly for $x_j \in I_j$. Note that up to terms of order $o_p(n^{-2/5})$ the second term on the right-hand side of (4.7) is asymptotically equal to a deterministic sequence. In the statement of Theorem 3 this sequence was called $\gamma_{n,j}$. The statement of Theorem 2.1 easily follows from (4.5)–(4.7).

We remark that in Assumption (A2) in [9] it was assumed that $\widetilde{p}_{jk}(x_j, x_k) = \widetilde{p}_{kj}(x_k, x_j)$ (in the notation of the current paper). Our choice of $\widetilde{p}_{jk}$ does not satisfy this symmetry constraint. It can be checked that Theorems 2 and 3 in [9] continue to hold when this symmetry constraint is dropped. Let us also mention that in their Assumption (A9) of Theorem 3 $\int \alpha_{n,j}(u) \widehat{p}_j(u) \, du = \gamma_{n,j} - o_p(\Delta_n)$ should be replaced by the correct assumption $[\int \widehat{p}_j(u) \, du]^{-1} \times \int \alpha_{n,j}(u) \widehat{p}_j(u) \, du = \gamma_{n,j} + o_p(\Delta_n)$.

It remains to show (4.3).

*Proof of* (4.3) *for* $S = B$. We first note that the following expansions hold:

$$(4.8) \quad \begin{aligned} \frac{\widetilde{p}_{jk}(x_j, x_k)}{\widetilde{p}_j(x_j)} &= \frac{\widehat{p}_{jk}(x_j, x_k)}{\widehat{p}_j(x_j)} + C_{K,j,2}(x_j) h_j^2 \frac{p_{jk}(x_j, x_k)}{p_j(x_j)^3} p_j'(x_j)^2 \\ &\quad - C_{K,j,2}(x_j) h_j^2 [D_{x_j} p_{jk}(x_j, x_k)] \frac{p_j'(x_j)}{p_j(x_j)^2} + o_p(n^{-2/5}), \end{aligned}$$

uniformly for $x_j \in I_j$ and $x_k \in I_k^-$, and

$$(4.9) \quad \frac{\widetilde{p}_{jk}(x_j, x_k)}{\widetilde{p}_j(x_j)} = \frac{\widehat{p}_{jk}(x_j, x_k)}{\widehat{p}_j(x_j)} + O_p(n^{-2/5}),$$

uniformly for $x_j \in I_j$ and $x_k \in I_k \setminus I_k^-$. These claims immediately follow from

$$(4.10) \quad \widehat{p}_j^*(x_j) = C_{K,j,2}(x_j) h_j^2 p_j'(x_j) + o_p(n^{-2/5}),$$

$$(4.11) \quad \widehat{p}_j^{**}(x_j) = C_{K,j,2}(x_j) h_j^2 p_j(x_j) + o_p(n^{-2/5}),$$

$$(4.12) \quad \widehat{p}_{jk}^*(x_j, x_k) = C_{K,j,2}(x_j) h_j^2 D_{x_j} p_{jk}(x_j, x_k) + o_p(n^{-2/5}),$$

uniformly for $x_j \in I_j$ and $x_k \in I_k^-$, and $\widehat{p}_j^*(x_j)$, $\widehat{p}_j^{**}(x_j)$, $\widehat{p}_{jk}^*(x_j, x_k)$ are all $O_p(n^{-2/5})$ uniformly for $x_j \in I_j$ and $x_k \in I_k \setminus I_k^-$.

Now, it follows that uniformly for $x_j \in I_j$

$$\int_{I_k} m_k(x_k) \widehat{p}_{jk}(x_j, x_k) \, dx_k$$

$$= \frac{1}{n}\sum_{i=1}^{n} K_{h_j}(x_j, X_j^i) \int_{I_k} m_k(x_k) L_{h_k}(x_k, X_k^i)\, dx_k$$

$$= \frac{1}{n}\sum_{i=1}^{n} K_{h_j}(x_j, X_j^i) m_k(X_k^i)$$

$$\quad - \frac{1}{2} C_K h_k^2 \left[\frac{1}{n}\sum_{i=1}^{n} K_{h_j}(x_j, X_j^i) m_k''(X_k^i)\right] + o_p(n^{-2/5})$$

$$= \frac{1}{n}\sum_{i=1}^{n} K_{h_j}(x_j, X_j^i) \int_{I_k} \frac{p_{jk}(X_j^i, x_k)}{p_j(X_j^i)} m_k(x_k)\, dx_k$$

$$\quad + \frac{1}{n}\sum_{i=1}^{n} K_{h_j}(x_j, X_j^i) \eta_{kj}^i$$

$$\quad - \frac{1}{2} C_K h_k^2 \int_{I_k} p_{jk}(x_j, x_k) m_k''(x_k)\, dx_k + o_p(n^{-2/5})$$

$$(4.13) \qquad = \widehat{p}_j(x_j) \int_{I_k} \frac{p_{jk}(x_j, x_k)}{p_j(x_j)} m_k(x_k)\, dx_k$$

$$\quad + C_{K,j,2}(x_j) h_j^2 \widehat{p}_j(x_j) \left[\left(D_{x_j} \int_{I_k} \frac{p_{jk}(x_j, x_k)}{p_j(x_j)} m_k(x_k)\, dx_k\right) \frac{p_j'(x_j)}{p_j(x_j)}\right.$$

$$\left. \quad + \frac{1}{2} D_{x_j}^2 \int_{I_k} \frac{p_{jk}(x_j, x_k)}{p_j(x_j)} m_k(x_k)\, dx_k\right]$$

$$\quad - \frac{1}{2} C_K h_k^2 \int_{I_k} p_{jk}(x_j, x_k) m_k''(x_k)\, dx_k + \frac{1}{n}\sum_{i=1}^{n} K_{h_j}(x_j, X_j^i)\eta_{kj}^i$$

$$\quad + o_p(n^{-2/5}).$$

Furthermore,

$$(4.14) \quad \begin{aligned} &\int_{I_k} \left[\frac{1}{2} C_K h_k^2 m_k''(x_k)\right] \frac{\widetilde{p}_{jk}(x_j, x_k)}{\widetilde{p}_j(x_j)}\, dx_k \\ &\qquad = \frac{1}{2} C_K h_k^2 \int_{I_k} \frac{p_{jk}(x_j, x_k)}{p_j(x_j)} m_k''(x_k)\, dx_k + o_p(n^{-2/5}). \end{aligned}$$

Using (4.8), (4.13) and (4.14) gives

$$\int_{I_k} \overline{m}_k^B(x_k) \frac{\widetilde{p}_{jk}(x_j, x_k)}{\widetilde{p}_j(x_j)}\, dx_k$$

$$(4.15) \qquad = \int_{I_k} \frac{p_{jk}(x_j, x_k)}{p_j(x_j)} m_k(x_k)\, dx_k + \frac{1}{n}\sum_{i=1}^{n} \frac{K_{h_j}(x_j, X_j^i)}{\widehat{p}_j(x_j)} \eta_{kj}^i$$

$$+ \frac{1}{2} C_{K,j,2}(x_j) h_j^2 \int_{I_k} \left[ D_{x_j}^2 \frac{p_{jk}(x_j, x_k)}{p_j(x_j)} \right] m_k(x_k) \, dx_k + o_p(n^{-2/5}).$$

Plugging (4.15) and (4.1) into the right-hand side of (4.3) gives uniformly for $x_j \in I_j$,

$$\check{m}_j^B(x_j) - \sum_{k \neq j} \int_{I_k} \overline{m}_k^B(x_k) \frac{\widetilde{p}_{jk}(x_j, x_k)}{\widetilde{p}_j(x_j)} \, dx_k$$

$$= m_j(x_j) + \frac{1}{2} C_{K,j,2}(x_j) h_j^2 m_j''(x_j) + o_p(n^{-2/5})$$

$$= \overline{m}_j^B(x_j) + o_p(n^{-2/5}).$$

This shows (4.3) for $S = B$.

*Proof of* (4.3) *for* $S = A$. We have to show for $k \neq j$ and $x_j \in I_j$,

$$\int_{I_k} \check{m}_k^A(x_k) \frac{\widetilde{p}_{jk}(x_j, x_k)}{\widetilde{p}_j(x_j)} \, dx_k = o_p(n^{-2/5}).$$

For this claim, it suffices to show that for $k \neq j$ and $x_j \in I_j$,

$$(4.16) \qquad \int_{I_k} \check{m}_k^A(x_k) \frac{\widehat{p}_{jk}(x_j, x_k)}{\widehat{p}_j(x_j)} \, dx_k = o_p(n^{-2/5}),$$

$$(4.17) \qquad \int_{I_k} \check{m}_k^A(x_k) \frac{p_{jk}(x_j, x_k)}{p_j(x_j)^3} p_j'(x_j)^2 \, dx_k = o_p(1),$$

$$(4.18) \qquad \int_{I_k} \check{m}_k^A(x_k) [D_{x_j} p_{jk}(x_j, x_k)] \frac{p_j'(x_j)}{p_j(x_j)^2} \, dx_k = o_p(1).$$

For the proof of (4.16)–(4.18), note that the left-hand sides of these equations can be written as $n^{-1} \sum_{i=1}^n w_i(x_j) \varepsilon_i$ where the weights $w_i(x_j)$ depend on $n, X^1, \ldots, X^n, x_j$, but not on $\varepsilon_1, \ldots, \varepsilon_n$. By standard smoothing theory it can be shown that in all three cases

$$\sup_{1 \leq i \leq n, x_j \in I_j} |w_i(x_j)| = O_p(1), \qquad \sup_{1 \leq i \leq n, x_j \in I_j} |w_i'(x_j)| = O_p(1).$$

These bounds imply

$$(4.19) \qquad \sup_{x_j \in I_j} \left| \frac{1}{n} \sum_{i=1}^n w_i(x_j) \varepsilon_i \right| = o_p(n^{-2/5}).$$

We give a short outline of the proof for (4.19).

Choose $C > 0$ and consider the event $E$ that $|w_i(x_j)| \leq C$ and $|w_i'(x_j)| \leq C$ for $1 \leq i \leq n$ and $x_j \in I_j$. We define

$$\overline{w}_i(x_j) = \begin{cases} w_i(x_j), & \text{on } E, \\ 1, & \text{elsewhere.} \end{cases}$$

Furthermore, for $\delta > 0$ small enough define

$$\overline{\varepsilon}_i = \varepsilon_i \mathbb{1}(|\varepsilon_i| \leq n^{1/2+\delta}) - E\varepsilon_i \mathbb{1}(|\varepsilon_i| \leq n^{1/2+\delta}).$$

Note that

$$P\left[\max_{1 \leq i \leq n} |\varepsilon_i| > n^{1/2+\delta}\right] \leq nP[|\varepsilon_1| > n^{1/2+\delta}]$$

$$\leq n^{-2\delta} E\varepsilon_1^2 = o(1)$$

and that

$$|E\varepsilon_1 \mathbb{1}(|\varepsilon_1| \leq n^{1/2+\delta})| = |E\varepsilon_1 \mathbb{1}(|\varepsilon_1| > n^{1/2+\delta})|$$

$$\leq n^{-1/2-\delta} E\varepsilon_1^2.$$

Therefore on $E$ we get

$$\frac{1}{n}\sum_{i=1}^{n} w_i(x_j)\varepsilon_i - \frac{1}{n}\sum_{i=1}^{n} w_i(x_j)\overline{\varepsilon}_i = O_p(n^{-1/2-\delta}) = o_p(n^{-2/5}).$$

Thus, it remains to show that

$$(4.20) \qquad \sup_{x_j \in I_j} \left|\frac{1}{n}\sum_{i=1}^{n} w_i(x_j)\overline{\varepsilon}_i\right| = o_p(n^{-2/5}).$$

For the proof of (4.20) we argue that

$$(4.21) \qquad \sup_{x_j \in I_j} \left|\frac{1}{n}\sum_{i=1}^{n} w'_i(x_j)\overline{\varepsilon}_i\right| = O_p(n^{1/2+\delta}),$$

and that for each $\Delta > 0$ there exist constants $C', C'' > 0$ such that

$$(4.22) \qquad \sup_{x_j \in I_j} P\left[\frac{1}{n}\sum_{i=1}^{n} w_i(x_j)\overline{\varepsilon}_i > \Delta n^{-2/5}\right] \leq C'\exp(-C''n^{1/10-\delta}).$$

We prove (4.22). On the event $E$,

$$P\left[\frac{1}{n}\sum_{i=1}^{n} w_i(x_j)\overline{\varepsilon}_i > \Delta n^{-2/5}\right]$$

$$\leq \exp(-n^{1/2-\delta}\Delta n^{-2/5})E\exp\left[n^{1/2-\delta}n^{-1}\sum_{i=1}^{n}\overline{w}_i(x_j)\overline{\varepsilon}_i\right]$$

$$\leq \exp(-\Delta n^{1/10-\delta})$$

$$\times \prod_{i=1}^{n}[1 + E(n^{-1-2\delta}\overline{w}_i(x_j)^2\overline{\varepsilon}_i^2 \exp(2n^{-1/2-\delta}|\overline{w}_i(x_j)|n^{1/2+\delta}))]$$

$$\leq \exp(-\Delta n^{1/10-\delta}) \prod_{i=1}^{n} [1 + C^2 \exp(2C) n^{-1-2\delta} E \overline{\varepsilon}_i^2]$$

$$\leq \exp(-\Delta n^{1/10-\delta}) \exp[n^{-2\delta} C^2 \exp(2C) E \overline{\varepsilon}_1^2]$$

$$= O(1) \exp(-\Delta n^{1/10-\delta}).$$

This shows (4.22) and completes the proof of Theorem 2.1.

**Acknowledgment.** We are grateful for the helpful and constructive comments of two reviewers.

## REFERENCES

[1] Buja, A., Hastie, T. and Tibshirani, R. (1989). Linear smoothers and additive models (with discussion). *Ann. Statist.* **17** 453–555. MR0994249

[2] Christopeit, N. and Hoderlein, S. (2003). Estimation of models with additive structure via local quasi-differencing. Preprint.

[3] Fan, J. and Gijbels, I. (1992). Variable bandwidth and local linear regression smoothers. *Ann. Statist.* **20** 2008–2036. MR1193323

[4] Fan, J., Härdle, W. and Mammen, E. (1998). Direct estimation of low-dimensional components in additive models. *Ann. Statist.* **26** 943–971. MR1635422

[5] Horowitz, J., Klemelä, J. and Mammen, E. (2006). Optimal estimation in additive regression models. *Bernoulli* **12** 271–298. MR2218556

[6] Horowitz, J. and Mammen, E. (2004). Nonparametric estimation of an additive model with a link function. *Ann. Statist.* **32** 2412–2443. MR2153990

[7] Linton, O. (1997). Efficient estimation of additive nonparametric regression models. *Biometrika* **84** 469–473. MR1467061

[8] Linton, O. and Nielsen, J. P. (1995). A kernel method of estimating structured nonparametric regression based on marginal integration. *Biometrika* **82** 93–100. MR1332841

[9] Mammen, E., Linton, O. and Nielsen, J. P. (1999). The existence and asymptotic properties of a backfitting projection algorithm under weak conditions. *Ann. Statist.* **27** 1443–1490. MR1742496

[10] Mammen, E., Marron, J. S., Turlach, B. and Wand, M. P. (2001). A general projection framework for constrained smoothing. *Statist. Sci.* **16** 232–248. MR1874153

[11] Mammen, E. and Park, B. U. (2005). Bandwidth selection for smooth backfitting in additive models. *Ann. Statist.* **33** 1260–1294. MR2195635

[12] Nielsen, J. P. and Sperlich, S. (2005). Smooth backfitting in practice. *J. R. Stat. Soc. Ser. B Stat. Methodol.* **67** 43–61. MR2136638

[13] Opsomer, J. D. (2000). Asymptotic properties of backfitting estimators. *J. Multivariate Anal.* **73** 166–179. MR1763322

[14] Opsomer, J. D. and Ruppert, D. (1997). Fitting a bivariate additive model by local polynomial regression. *Ann. Statist.* **25** 186–211. MR1429922

[15] Opsomer, J. D. and Ruppert, D. (1998). A fully automated bandwidth selection method for fitting additive models. *J. Amer. Statist. Assoc.* **93** 605–619. MR1631333

[16] Tjøstheim, D. and Auestad, B. H. (1994). Nonparametric identification of nonlinear time series: Projections. *J. Amer. Statist. Assoc.* **89** 1398–1409. MR1310230

Department of Economics
University of Mannheim L 7, 3-5
68131 Mannheim
Germany
E-mail: emammen@rumms.uni-mannheim.de

Department of Statistics
Seoul National University
Seoul 151-747
Korea
E-mail: bupark@stats.snu.ac.kr